\documentclass[12pt,reqno]{amsart}

\usepackage{times}
\usepackage{amssymb,amsmath,amsfonts,amsthm,mathrsfs}
\usepackage{hyperref}
\hypersetup{hidelinks}
\usepackage[inner=1.0in,outer=1.0in,bottom=1.0in,top=1.0in]{geometry}

\newtheorem{theorem}{Theorem}[section]

\newtheorem{proposition}[theorem]{Proposition}
\newtheorem{lemma}[theorem]{Lemma}
\newtheorem{corollary}[theorem]{Corollary}

\newtheorem{remark}[theorem]{Remark}
\numberwithin{equation}{section}

\newcommand{\eps}{\varepsilon}

\newcommand{\HH}{\mathbb H}

\newcommand{\RR}{\mathbb R}

\newcommand{\ZZ}{\mathbb Z}

\DeclareMathOperator{\vol}{vol}
\DeclareMathOperator{\sym}{sym}

\DeclareMathOperator{\PSL}{PSL}

\DeclareMathOperator{\supp}{supp}

\title{Quantitative nodal domain count for Maa\ss\ forms}

\author[P. Zenz]{Peter Zenz}
\address{Department of Mathematics, Brown University, Providence, RI 02912}
\email{peter\_zenz@brown.edu}

\begin{document}

\begin{abstract}
We consider a real even Hecke–Maaß cusp form $\phi$ on the modular surface, with Laplace eigenvalue $1/4+T^2$. We prove that its restriction to any fixed compact segment $\beta$ of the vertical geodesic above $i$ has $\gg_{\beta,\varepsilon} T^{1/6-\varepsilon}$ sign changes, yielding the same lower bound for the number of distinct inert nodal domains meeting the segment.
Building on ideas of Jang–Jung and Jung–Young, we study correlations between values at points on the vertical geodesic separated by distances of order $T^{-1}$. Using arithmetic quantum unique ergodicity and Dirichlet-polynomial estimates, we show that the correlation is negative for a suitable choice of separation. An elementary counting argument and the Iwaniec–Sarnak supremum bound then give the claimed lower bound.
\end{abstract}

\maketitle

\section{Introduction}
\label{sec:introduction}

Nodal domains describe how an eigenfunction divides a surface into
connected regions of constant sign. Understanding how their number grows
with the eigenvalue is a basic problem in spectral geometry.
We study sign changes of even Hecke--Maa\ss\ cusp forms on a prescribed
compact segment of the imaginary axis. By reflection symmetry, lower bounds for these sign changes give
lower bounds for the number of nodal domains meeting the segment.

We work on $X=\PSL_2(\ZZ)\backslash\HH$, with hyperbolic area
$d\mu=dx\,dy/y^2$ and $\vol(X)=\pi/3$. Throughout, $\phi$ is a real-valued
even Hecke--Maa\ss\ cusp form satisfying
\[
 \Delta_+\phi=(T^2+1/4)\phi,\qquad
 \Delta_+=-y^2(\partial_x^2+\partial_y^2),\qquad
 \int_X|\phi|^2\,d\mu=1,
\]
where $T>0$, and $\phi(-x+iy)=\phi(x+iy)$.
We fix $1<y_0<y_1<\infty$ and set
$\beta=\{iy:y_0\le y\le y_1\}$.
We write $K_\beta(\phi)$ for the number of distinct sign-changing
zeros of the restriction in the interior of $\beta$.
A nodal domain is a connected component of $X\setminus\{\phi=0\}$.
The reflection $\sigma(z)=-\overline z$ permutes these domains.
A domain $D$ is called inert if $\sigma(D)=D$; otherwise $D$ and
$\sigma(D)$ form a pair of split domains. Every domain meeting
$\beta$ is inert, since it contains a point fixed by $\sigma$.
We write $N_{\mathrm{in}}^\beta(\phi)$ for the number of distinct
such domains \cite{GRS13}.

\begin{theorem}\label{thm:main}
For every fixed segment $\beta$ as above and every $\eps>0$,
every form $\phi$ as above with sufficiently large spectral
parameter $T$ satisfies
\begin{equation}\label{eq:main-counts}
 K_\beta(\phi)\gg_{\beta,\eps}T^{1/6-\eps},\qquad
 N_{\mathrm{in}}^\beta(\phi)\gg_{\beta,\eps}T^{1/6-\eps}.
\end{equation}
\end{theorem}

The random-wave model suggests that the number of inert nodal
domains meeting a fixed compact segment $\beta$ should be of
order $T$ \cite{BS02,GRS13}. For the total number of nodal domains
on $X$, the Bogomolny--Schmit prediction is of order $T^2$
\cite{BS02}.

Ghosh, Reznikov and Sarnak proved a lower bound
$\gg T^{1/12-\eps}$ for inert nodal domains meeting a sufficiently
long fixed compact vertical segment, assuming the Lindel\"of hypothesis
for $L(s,\phi)$ \cite[Thm.~1.8]{GRS13}.
Jung obtained an unconditional lower bound of order $T^{1/8-\eps}$
for nodal domains meeting any prescribed fixed compact vertical segment,
for almost all even forms in short spectral windows \cite{JungSC}.
In the compact setting, Jang and Jung proved unconditionally that
the nodal-domain count tends to infinity along every sequence of
Hecke--Maa\ss\ eigenforms on an arithmetic triangle surface whose
eigenvalues tend to infinity \cite{JJ18}. Their argument uses quantum
unique ergodicity (QUE) \cite{Lin06} and a separating reflection symmetry.
Theorem~\ref{thm:main} gives the lower bound $T^{1/6-\eps}$ for both
counts on every fixed segment $\beta$, for every even Hecke--Maa\ss\
cusp form with sufficiently large spectral parameter.

Under Lindel\"of and a conjectural fourth-moment bound on the vertical
geodesic, Ki obtains the essentially optimal lower bound $T^{1-\eps}$
on sufficiently long fixed compact segments of the imaginary axis
\cite{Ki26}.
Kelmer, Kontorovich and Lutsko obtain the lower bound $T^{1-\eps}$
for sign changes on any fixed compact segment of the imaginary axis,
assuming $\int_\beta|\phi(iy)|^p\,dy/y\ll_{\beta,p,\eps}T^\eps$
for some fixed $2<p<\infty$ and every $\eps>0$, uniformly in the family under
consideration, together with a mean-square conjecture for $L(s,\phi)$
\cite{KKL25}. Corollary~\ref{cor:restricted-moment-count} gives the
same lower bound for both counts under the restricted-moment assumption
alone.

Our main analytic result evaluates a difference of correlations on the
vertical geodesic. We fix a nonzero
$\psi\in C_c^\infty((y_0,y_1))$ with $0\le\psi\le1$. Following Jung and Young \cite{JY19}, we define, for fixed $\alpha\in\RR$ and
sufficiently large $T$,
\begin{equation}\label{eq:arithmetic-correlation}
 I_\phi(\alpha)=\int_0^\infty
 \psi(y)\psi((1+\alpha/T)y)
 \phi(iy)\phi(i(1+\alpha/T)y)\,\frac{dy}{y}.
\end{equation}
We put $H_\psi=(3/\pi)\int_0^\infty\psi(y)^2\,dy/y$, so that $H_\psi$ is the cutoff integral divided by $\vol(X)$. The unshifted value $I_\phi(0)$ is
the weighted restriction mass. We evaluate the
difference $I_\phi(\alpha)-\cos\alpha\,I_\phi(0)$.

\begin{theorem}\label{prop:correlation}
For every fixed $\alpha\in\RR$ and every fixed $\psi$ as above,
\begin{equation}\label{eq:arithmetic-defect}
 I_\phi(\alpha)-\cos\alpha\,I_\phi(0)
 =2H_\psi\bigl(J_0(\alpha)-\cos\alpha\bigr)
       +o_{\alpha,\psi}(1)
\end{equation}
as $T\to\infty$, uniformly over the forms $\phi$ specified above.
Here $J_0$ denotes the Bessel function of order zero.
\end{theorem}
The $o(1)$ remainder comes from the qualitative QUE input in the bulk
calculation. A suitable quantitative QUE estimate would give a quantitative
error term here.
At the shift
\begin{equation}\label{eq:chosen-shift}
 \alpha_*=\frac{3\pi}{2},
\end{equation}
the coefficient of $I_\phi(0)$ vanishes and $J_0(\alpha_*)=-0.265857\ldots<-1/4$.
Consequently, $I_\phi(\alpha_*)$ is bounded above by a fixed negative constant for all forms with sufficiently large $T$. 
This forces many sign changes, since each zero can account for only a small portion of the negative mass.

More precisely, we put $S_\beta(\phi)=\sup_\beta|\phi|$. A negative contribution to
$I_\phi(\alpha_*)$ requires a sign-changing zero between $y$ and
$(1+\alpha_*/T)y$. For any given zero, this can occur only for $y$
in an interval of length $O_\beta(T^{-1})$. Since the product in the
correlation has absolute value at most $S_\beta(\phi)^2$, the total
negative contribution has magnitude at most
$O_\beta(K_\beta(\phi)S_\beta(\phi)^2/T)$. The correlation formula
gives a fixed positive lower bound for this magnitude, so
\[
 K_\beta(\phi)\gg_\psi\frac{T}{S_\beta(\phi)^2}.
\]
The Iwaniec--Sarnak supremum bound therefore gives $T^{1/6-\eps}$
sign changes, and reflection symmetry gives the same lower bound
for the number of nodal domains \cite[(6)]{GRS13}.
Lemma~\ref{lem:negative-pairs} makes the counting argument precise,
and we give the complete deduction in Section~\ref{sec:consequences}.

Our argument draws on Jung and Young's correlation method and Jang and
Jung's use of QUE to detect sign changes. The cosine--Bessel combination
in \eqref{eq:arithmetic-defect} already appears in Jang and Jung's
characteristic-function argument \cite{JJ18}. Their method is closely
related to the weighted Cauchy-data restriction theorem of Christianson,
Toth and Zelditch \cite{CTZ}. We explain this relation and the common
cancellation near $t=\pm T$ in Section~\ref{sec:outline}.

Jung and Young \cite{JY19} obtained quantitative nodal-domain bounds
for Eisenstein series from a correlation asymptotic.
Their correlation, divided by its unshifted mass, tends to
$J_0(\alpha)$ for each fixed $\alpha$. Their argument uses the
explicit Fourier coefficients of Eisenstein series. For general
Hecke--Maa\ss\ cusp forms, arithmetic QUE does not directly give
an asymptotic for $I_\phi(0)$
\cite[\S1.4, \S1.6]{YoungGeo}. We instead evaluate the difference
in \eqref{eq:arithmetic-defect}, without requiring an asymptotic
for this mass. Section~\ref{sec:outline} explains the difficulty
and the cancellation that makes this possible.

Our counting lemma also gives the lower bound $T^{1/3-\eps}$
for both counts in the Eisenstein case on every fixed $\beta$
(see Corollary~\ref{cor:eisenstein-count}).

\subsection*{Acknowledgments}
We carried out the initial part of this work while holding a Visiting
Assistant Professorship at Stanford University. We thank
the Department of Mathematics for its hospitality. We also thank Kannan Soundararajan, Junehyuk Jung and Matt Young for helpful discussions and interest.

\subsection*{Use of AI tools}
Discussions with ChatGPT (OpenAI) contributed to the development of the counting argument in Lemma~\ref{lem:negative-pairs} and the
proof in Appendix~\ref{app:coefficient-que}. Generative AI also assisted with the exposition and editing of the manuscript. We take responsibility for the mathematical content.

\subsection*{Notation and normalization}
We write $e(x)=e^{2\pi ix}$. We indicate the dependence of implied constants
by subscripts or state it in the text. We use the Fourier expansion of $\phi$ (see \cite[\S3, (13)--(16)]{GRS13}):
\begin{equation}
 \phi(x+iy)=\rho_\phi(1)\sqrt y
 \sum_{n\ne0}\lambda_\phi(|n|)K_{iT}(2\pi|n|y)e(nx),
 \qquad
 |\rho_\phi(1)|^2=\frac{2\cosh(\pi T)}{L(1,\sym^2\phi)}.
 \label{eq:fourier-normalization}
\end{equation}
Here $\rho_\phi(1)$ is real, $\lambda_\phi(n)$ are the Hecke eigenvalues
with $\lambda_\phi(1)=1$, and $K_{iT}$ denotes the $K$-Bessel function.
The identity for $|\rho_\phi(1)|^2$ follows from the $L^2$-normalization
$\int_X|\phi|^2\,d\mu=1$ by Rankin--Selberg unfolding
\cite[\S2, (2.1), (2.16), (2.19), Lem.~2.2]{YoungThin}.

\section{Outline}
\label{sec:outline}

For $f(y)=\psi(y)\phi(iy)$, we put
$\mathcal M(it)=\int_0^\infty f(y)y^{it}\,dy/y$.
We apply the Mellin--Parseval formula to express the correlation
difference in \eqref{eq:arithmetic-defect} as
\[
 I_\phi(\alpha)-\cos\alpha\,I_\phi(0)
 =\frac1{2\pi}\int_\RR
 \bigl(\cos(t\log(1+\alpha/T))-\cos\alpha\bigr)
 |\mathcal M(it)|^2\,dt.
\]
For $|t|\le T$, the factor in parentheses is essentially
$\cos(\alpha t/T)-\cos\alpha$. This multiplicative
factor reduces the contribution from parameters $t$ that are close to $\pm T$.
Indeed, for $t$ such that $H\le T-|t|\le2H$, with $1\leq H=o(T)$,
we gain an additional factor of size $H/T$.

To study the correlation, we approximate $\mathcal M(it)$ by
an explicit gamma factor times a Dirichlet polynomial:
\[
 \mathcal M(it)\approx G_T(it)P_T(t),\qquad
 P_T(t)=\sum_{n\ge1}\frac{\lambda_\phi(n)}{n^{1/2+it}}
 \psi\left(\frac{\sqrt{T^2-t^2}}{2\pi n}\right),
\]
with $|G_T(it)|^2\sim
2\pi/(L(1,\sym^2\phi)\sqrt{T^2-t^2})$ for
$|t|\le T-T^\delta$, for any fixed $\delta>0$.
We choose $0<\delta<1/6$. The correlation estimate is then morally equivalent to the following
mean value problem of Dirichlet polynomials
\[
 \frac1{L(1,\sym^2\phi)}
 \int_{|t|\le T-T^\delta}
 \frac{\big(\cos(\alpha t/T)-\cos\alpha\big)}{\sqrt{T^2-t^2}}
 \Bigg|\sum_{n \sim \sqrt{T^2-t^2}}\frac{\lambda_\phi(n)}{n^{1/2+it}}\Bigg|^2\,dt.
\]
We divide the integration range into a bulk region and boundary blocks approaching
$t=\pm T$. We treat the remaining neighborhoods
$\bigl||t|-T\bigr|\le T^\delta$ and the exterior
$|t|>T+T^\delta$ separately below.

We fix a small $\eta>0$. In the bulk $|t|\le(1-\eta)T$,
the Dirichlet polynomial has length comparable to $T$.
Expanding its squared absolute value gives smoothly weighted sums
of $\lambda_\phi(n)\lambda_\phi(n+h)$ with $n\asymp T$.
The diagonal gives the main term, while full phase-space QUE
gives cancellation for each fixed nonzero shift $h$.
The oscillation in $t$ gives a summable bound for large
$|h|$, so qualitative convergence for fixed shifts suffices.
We carry out this calculation for the actual polynomial and its moving
cutoff in Section~\ref{sec:bulk}.

For the boundary frequencies, we divide
$T^\delta<T-|t|\le\eta T$ into dyadic blocks of scale $H$.
The lower cutoff separates off the intervals
$\bigl||t|-T\bigr|\le T^\delta$ around $t=\pm T$,
and the upper cutoff separates the boundary blocks from the bulk.
First, we keep $\eta$ fixed in the QUE limit so that the coefficient weights
are fixed, and later we let $\eta\downarrow0$.
To see why the factor $\cos(\alpha t/T)-\cos\alpha$ is needed,
we first consider a boundary block of the unshifted mass, where this
factor is absent. We suppress fixed constants and smooth weights
to obtain an expression of the form
\[
 \frac1{L(1,\sym^2\phi)}
 \int_{T-2H}^{T-H}\frac1{\sqrt{TH}}
 \left|\sum_{n\sim\sqrt{TH}}
       \frac{\lambda_\phi(n)}{n^{1/2+it}}\right|^2dt.
\]
The polynomial length $\sqrt{TH}$ exceeds the interval length $H$.
The standard Montgomery--Vaughan mean-value theorem \cite{MV74} bounds this expression by
a constant multiple of
\[
 \frac{H+\sqrt{TH}}{\sqrt{TH}}\,
 \frac1{L(1,\sym^2\phi)}
 \sum_{n\sim\sqrt{TH}}\frac{|\lambda_\phi(n)|^2}{n}.
\]
The normalized coefficient average is expected to have size one.
At this expected size, the mean-value estimate gives only $O(1)$
per block, and summing over the dyadic blocks gives
$O(\log T)$. Overcoming this obstruction is the main bottleneck for proving an asymptotic
formula for $I_\phi(0)$. We circumvent this problem by considering $I_\phi(\alpha)-\cos\alpha\,I_\phi(0)$. Here the additional
factor $H/T$ gives the estimate
\begin{equation}\label{dyadic_summable}
 \sum_{\substack{H\ \mathrm{dyadic}\\T^\delta\lesssim H\lesssim\eta T}}
 \frac HT\left(1+\sqrt{\frac HT}\right)\ll\eta.
\end{equation}
In Section~\ref{sec:arithmetic-boundary}, we prove a slightly weaker but sufficient estimate. On the intervals $\bigl||t|-T\bigr|\le T^\delta$,
the small multiplicative factor in the correlation difference
and the compact supremum bound give a negligible
contribution because $\delta<1/6$. Gamma-function decay controls the
exterior. In Section~\ref{sec:arithmetic-completion}, we combine these estimates
with the bulk limit to prove Theorem~\ref{prop:correlation}.

\subsection*{Comparison with earlier methods}
On compact manifolds with ergodic geodesic flow, Toth and Zelditch
proved equidistribution of restrictions along a density-one subsequence
under a microlocal asymmetry condition \cite{TZ13}. For a reflection
geodesic, the relevant comparison is the weighted Cauchy-data theorem
of Christianson, Toth and Zelditch \cite{CTZ}. It derives
equidistribution of a weighted sum of the quadratic forms of the
restriction and its normalized normal derivative from equidistribution
of the ambient eigenfunctions. The weight on the restriction term
vanishes at unit tangential directions. This formulation accommodates
reflection symmetry and is closely related to Jang and Jung's use of
QUE to detect sign changes.

Jang and Jung use a Rellich identity and QUE to evaluate differences
between the squared $L^2$ norm of a localized restriction and those
of its normalized tangential derivatives \cite{JJ18}.
For the first derivative, the corresponding expression in our
variables satisfies the exact Mellin--Parseval identity
\begin{equation}
 \begin{aligned}
 &\int_0^\infty\left((\psi(y)\phi(iy))^2
 -\frac{|y\partial_y(\psi(y)\phi(iy))|^2}{T^2+1/4}\right)\frac{dy}{y}\\
 &\qquad=\frac1{2\pi}\int_\RR
 \left(1-\frac{t^2}{T^2+1/4}\right)|\mathcal M(it)|^2\,dt.
 \end{aligned}
 \label{eq:mellin-derivative}
\end{equation}
Its leading factor $1-t^2/T^2$ contributes $H/T$ on a boundary
block. The endpoints $t=\pm T$ correspond to the normalized tangential
frequencies $\xi=\pm1$, where $\xi=t/T$. We can compare this with
the Cauchy-data weight $1-|\xi|^2$ of Christianson, Toth and Zelditch,
which also vanishes there \cite[Thm.~1, \S2]{CTZ}. In the compact setting,
Jang and Jung evaluate these differences using QUE, then use convergence
of moments and Bochner's theorem to detect sign changes. For our correlation, we notice that
the factor $\cos(\alpha t/T)$ for $\alpha=3\pi/2$ contributes a similar saving of size $H/T$.

Jang and Jung's moment argument already contains the cosine--Bessel
combination in \eqref{eq:arithmetic-defect} \cite{JJ18}.
After normalizing the restriction mass and passing to a subsequence,
their limiting frequency measure $\mu_a$ satisfies
\[
 \int_\RR e^{i\alpha\xi}\,d\mu_a(\xi)-\cos\alpha
 =a\bigl(J_0(\alpha)-\cos\alpha\bigr),\qquad 0\le a\le1.
\]
In our normalization, $a$ corresponds to a subsequential limit of
$2H_\psi/I_\phi(0)$. The equal masses at $\xi=1$ and $\xi=-1$ contribute
zero to this difference. Formula~\eqref{eq:arithmetic-defect} evaluates
the unnormalized correlation difference with an absolute $o(1)$ remainder.

At the shifts $\alpha=\pi/2+k\pi$, with fixed $k\in\ZZ$, the multiplicative factor
$\cos(\alpha t/T)$ itself vanishes at $t=\pm T$. We briefly explain from a different perspective why
these values are special and how Bessel functions enter. To do so, we temporarily replace the moving
coefficients of $P_T(t)$ by real coefficients $b_n$ independent of $t$. For the resulting Dirichlet
polynomial $P(t)=\sum_n b_n n^{-it}$, we open the square and pair the
terms to obtain $|P(t)|^2=\sum_{m,n}b_mb_n\cos(u_{m,n}t/T)$,
where $u_{m,n}=T\log(m/n)$. We define
\begin{align*}
 K_\alpha(u)
 &:=\frac1\pi\int_{-T}^T
 \frac{\cos(\alpha t/T)\cos(ut/T)}{\sqrt{T^2-t^2}}\,dt\\
 &=\frac12\bigl(J_0(u+\alpha)+J_0(u-\alpha)\bigr).
\end{align*}
Thus
\[
 \frac1\pi\int_{-T}^T
 \frac{\cos(\alpha t/T)|P(t)|^2}{\sqrt{T^2-t^2}}\,dt
 =\sum_{m,n}b_mb_nK_\alpha(u_{m,n}).
\]
For comparison, we insert the factor $1-t^2/T^2$ and obtain
\[
 \frac1\pi\int_{-T}^T
 \left(1-\frac{t^2}{T^2}\right)
 \frac{\cos(ut/T)}{\sqrt{T^2-t^2}}\,dt
 =\frac{J_1(u)}u=\frac{J_0(u)+J_2(u)}2,
\]
with value $1/2$ at $u=0$. These identities follow from
\cite[10.9.4, 10.6.1]{DLMF}. For fixed $\alpha$, the asymptotic
expansion \cite[10.17.3]{DLMF} gives
\[
 K_\alpha(u)=\left(\frac2{\pi u}\right)^{1/2}
 \cos\left(u-\frac\pi4\right)\cos\alpha
 +O_\alpha(u^{-3/2})\qquad(u\to+\infty).
\]
Hence $K_\alpha(u)-\cos\alpha\,K_0(u)=O_\alpha(|u|^{-3/2})$,
the same decay as $J_1(u)/u$. When $\cos\alpha=0$, the leading
oscillations of the two shifted Bessel terms cancel, gaining a factor
$|u|^{-1}$. Thus the kernel for the correlation itself has the improved
decay at these shifts. In the proof with moving coefficients, the same
cancellation makes the boundary estimates summable (see \eqref{dyadic_summable}). 
At $\alpha=\alpha_*=3\pi/2$, the diagonal value
$K_\alpha(0)=J_0(\alpha)$ is also negative. This choice therefore
gives both the cancellation and the negative main term needed
for our sign-change argument.

\section{From negative correlation to nodal domains}
\label{sec:consequences}
\label{sec:counting}

We derive the counting consequences of Theorem~\ref{prop:correlation},
whose proof is completed in Section~\ref{sec:arithmetic-completion}.
To turn a negative correlation into a count of sign changes, we bound
how much each zero can contribute. If $u(s)$ and $u(s+\delta)$ have
opposite signs, then a sign-changing zero lies between $s$ and
$s+\delta$. A given zero $z$ lies between these points only when
$z-\delta<s<z$, so the corresponding values of $s$ form an interval
of length $\delta$. If $|u|\le S$, the negative contribution from
this interval has magnitude at most $\delta S^2$. Summing over the
sign-changing zeros gives the following elementary estimate.

\begin{lemma}\label{lem:negative-pairs}
Let $u$ be real-valued and analytic on a connected neighborhood of a
compact interval $I$, and suppose that $u$ is not identically zero.
Let $K$ be the number of distinct sign-changing zeros of $u$ in
$I^\circ$ and let $S=\sup_I|u|$. Let $w$ be measurable, with $0\le w\le1$ and support in $I^\circ$,
and extend $v=wu$ by zero outside $I$. Then, for every $\delta>0$,
\begin{equation}\label{eq:negative-pair-bound}
 \int_{\RR}v(s)v(s+\delta)\,ds\ge-\delta KS^2.
\end{equation}
In particular, if the integral is at most $-D<0$, then
$K\ge D/(\delta S^2)$.
\end{lemma}

\begin{proof}
Since $u$ is real analytic and not identically zero, it has only
finitely many zeros in $I$. Put
$E_\delta=\{s:v(s)v(s+\delta)<0\}$. If $s\in E_\delta$, then
both endpoints lie in $I^\circ$, the cutoff is positive there, and
$u(s)$ and $u(s+\delta)$ have opposite signs. Hence some sign-changing
zero $z_j$ lies in $(s,s+\delta)$, so
\[
 E_\delta\subseteq\bigcup_{j=1}^K(z_j-\delta,z_j),
 \qquad |E_\delta|\le K\delta.
\]
The integrand is nonnegative outside $E_\delta$ and is at least
$-S^2$ on $E_\delta$. This proves \eqref{eq:negative-pair-bound}.
\end{proof}

For related bounded-variation translation estimates, see Galerne and
Lachi\`eze-Rey \cite{GLR15}. Short-integral criteria for sign changes
appear in Littlewood \cite{Lit66}, Ki \cite{Ki26}, and Kelmer,
Kontorovich and Lutsko \cite{KKL25}.

\begin{proof}[Proof of Theorem~\ref{thm:main}]
Fix $\psi$ as in Theorem~\ref{prop:correlation} for the prescribed
segment $\beta$. At $\alpha_*=3\pi/2$, the coefficient of $I_\phi(0)$ vanishes
and $J_0(\alpha_*)<-1/4$. Thus
\[
 I_\phi(\alpha_*)=2H_\psi J_0(\alpha_*)+o_\psi(1),
\]
and for every sufficiently large form,
\begin{equation}\label{eq:negative-correlation}
 I_\phi(\alpha_*)\le-\frac{H_\psi}{4}.
\end{equation}
Here the threshold depends only on the fixed cutoff: convergence in
Theorem~\ref{prop:correlation} holds along every sequence of forms
with $T\to\infty$.

Set
\[
 I=[\log y_0,\log y_1],\qquad u(s)=\phi(ie^s),\qquad
 w(s)=\psi(e^s),\qquad \Delta_T=\log(1+\alpha_*/T).
\]
The change of variables $y=e^s$ gives
$I_\phi(\alpha_*)=\int_{\RR}(wu)(s)(wu)(s+\Delta_T)\,ds$.
The restriction $u$ is real analytic, and
\eqref{eq:negative-correlation} shows that it is not identically zero.
Lemma~\ref{lem:negative-pairs} therefore gives
\begin{equation}\label{eq:quantitative-count}
 K_\beta(\phi)\ge
 \frac{H_\psi}{4\log(1+\alpha_*/T)S_\beta(\phi)^2},
 \qquad S_\beta(\phi)=\sup_\beta|\phi|.
\end{equation}
Since $\log(1+\alpha_*/T)\asymp T^{-1}$, the Iwaniec--Sarnak supremum
bound \cite{IS95}, recalled in \eqref{eq:supnorm}, with exponent
$5/12+\eps/2$, yields
\[
 K_\beta(\phi)\gg_{\beta,\eps}T^{1/6-\eps}.
\]

We apply the topological argument of Ghosh, Reznikov and Sarnak
\cite[(6)]{GRS13} to obtain
\begin{equation}\label{eq:topology}
 N_{\mathrm{in}}^\beta(\phi)\ge1+\tfrac12K_\beta(\phi).
\end{equation}
This proves the second assertion.
\end{proof}

We can replace the supremum bound by a restricted $L^p$ moment.
\begin{corollary}\label{cor:restricted-moment-count}\label{cor:fourth-moment-count}
For every fixed segment $\beta$ as above and every fixed $2<p<\infty$,
every normalized cusp form $\phi$ as above with sufficiently large
spectral parameter $T$ satisfies
\begin{equation}\label{eq:restricted-lp-count}
 K_\beta(\phi)\gg_{\beta,p}
 \frac{T}{\left(\displaystyle\int_\beta|\phi(iy)|^p\,\frac{dy}{y}
          \right)^{2/(p-2)}}.
\end{equation}
If, for one fixed $2<p<\infty$, a family of these forms satisfies
\[
 \int_\beta|\phi(iy)|^p\,\frac{dy}{y}
 \ll_{\beta,p,\eps}T^\eps
\]
for every $\eps>0$ uniformly throughout the family, then for every
$\eps>0$ all sufficiently large members of the family satisfy
\[
 K_\beta(\phi)\gg_{\beta,p,\eps}T^{1-\eps},\qquad
 N_{\mathrm{in}}^\beta(\phi)\gg_{\beta,p,\eps}T^{1-\eps}.
\]
\end{corollary}

\begin{proof}
Fix $\psi$ for $\beta$ as above, put
$\Delta_T=\log(1+\alpha_*/T)$ and
$v(s)=\psi(e^s)\phi(ie^s)$, and extend $v$ by zero outside
$[\log y_0,\log y_1]$. For
$E=\{s:v(s)v(s+\Delta_T)<0\}$, the covering in the proof of
Lemma~\ref{lem:negative-pairs} gives
$|E|\le K_\beta(\phi)\Delta_T$.
By \eqref{eq:negative-correlation}, H\"older's inequality with exponents
$p/2$ and $p/(p-2)$, and Cauchy--Schwarz, we obtain
\[
\begin{aligned}
 \frac{H_\psi}{4}
 &\le \int_E|v(s)v(s+\Delta_T)|\,ds\\
 &\le |E|^{1-2/p}
       \left(\int_{\RR}|v(s)|^{p/2}|v(s+\Delta_T)|^{p/2}\,ds\right)^{2/p}\\
 &\le \bigl(K_\beta(\phi)\Delta_T\bigr)^{1-2/p}
       \left(\int_{\RR}|v(s)|^p\,ds\right)^{2/p}\\
 &\le \bigl(K_\beta(\phi)\Delta_T\bigr)^{1-2/p}
       \left(\int_\beta|\phi(iy)|^p\,\frac{dy}{y}\right)^{2/p}.
\end{aligned}
\]
The third line also uses translation invariance on $\RR$, and the
last uses $0\le\psi\le1$ and its support in
$(y_0,y_1)$. The integral in the last line is therefore positive.
Raising to the power $p/(p-2)$ and rearranging gives
\[
 K_\beta(\phi)\ge
 \frac{(H_\psi/4)^{p/(p-2)}}
 {\Delta_T\left(\displaystyle\int_\beta|\phi(iy)|^p\,\frac{dy}{y}
          \right)^{2/(p-2)}}.
\]
Since $\Delta_T\asymp T^{-1}$, this proves
\eqref{eq:restricted-lp-count}. For $p=4$, the estimate becomes
\begin{equation}\label{eq:fourth-moment-count}
 K_\beta(\phi)\ge
 \frac{H_\psi^2}
 {16\Delta_T\displaystyle\int_\beta|\phi(iy)|^4\,\frac{dy}{y}}
 \gg_\beta\frac{T}{\displaystyle\int_\beta|\phi(iy)|^4\,\frac{dy}{y}}.
\end{equation}
Finally, applying the moment hypothesis with exponent $\eps(p-2)/2$
gives the sign-count bound, and \eqref{eq:topology} gives the
inert-domain bound.
\end{proof}

For Eisenstein series, Jung and Young \cite{JY19} evaluated the
fixed-shift correlation needed here. Combined with
Lemma~\ref{lem:negative-pairs} and Blomer's compact supremum bound
\cite{Blomer20}, their formula gives $T^{1/3-\eps}$ sign changes
on every fixed segment $\beta$ for every sufficiently large $T$. This improves the exponent $1/51$ in their nodal domain counting result.
Write
\[
 E_T^*(z)=\frac{\theta(1/2+iT)}{|\theta(1/2+iT)|}E(z,1/2+iT),
 \qquad \theta(s)=\pi^{-s}\Gamma(s)\zeta(2s),
\]
where $E(z,s)$ is the standard Eisenstein series for $\PSL_2(\ZZ)$.
Then $E_T^*$ is real-valued and $|E_T^*|=|E(z,1/2+iT)|$.

\begin{corollary}\label{cor:eisenstein-count}\label{rem:eisenstein-count}
For every fixed segment $\beta$ as above and every $\eps>0$,
all sufficiently large real $T$ satisfy
\[
 K_\beta(E_T^*)\gg_{\beta,\eps}T^{1/3-\eps},\qquad
 N_{\mathrm{in}}^\beta(E_T^*)\gg_{\beta,\eps}T^{1/3-\eps}.
\]
\end{corollary}

\begin{proof}
Fix $\beta$ and $\psi$ as above, and define $I_{E_T^*}(\alpha)$ by
\eqref{eq:arithmetic-correlation}. The fixed-cutoff case of the
Jung--Young correlation formula \cite[Thm.~3.2]{JY19}, which allows
any fixed compact interval, gives
\[
 I_{E_T^*}(\alpha_*)
 =2H_\psi\log(T^2+1/4)J_0(\alpha_*)+o_\psi(\log T)
 \le-c_\psi\log T
\]
for every sufficiently large $T$, with $c_\psi>0$.
Lemma~\ref{lem:negative-pairs}, applied to $u(s)=E_T^*(ie^s)$,
$w(s)=\psi(e^s)$ and displacement $\Delta_T=\log(1+\alpha_*/T)$,
therefore gives
\[
 K_\beta(E_T^*)\gg_\psi
 \frac{T\log T}{\sup_\beta|E_T^*|^2}.
\]
The Epstein series for $Q_y(c,d)=yc^2+d^2/y$ equals
$2\zeta(2s)E(iy,s)$, by separating primitive lattice vectors.
Thus Blomer's bound \cite[Thm.~1, $k=2$]{Blomer20}, uniform for
$y$ in our fixed interval, together with the standard bound
$|\zeta(1+2iT)|^{-1}=T^{o(1)}$ \cite[(4.5)]{JY19}, gives
\[
 \sup_\beta|E_T^*|\ll_{\beta,\eps}T^{1/3+\eps/2},
\]
and hence $K_\beta(E_T^*)\gg_{\beta,\eps}T^{1/3-\eps}$.
The local reflection argument giving the nodal domain count \eqref{eq:topology} also applies
to Eisenstein series
\cite[Cor.~1.2]{JY19}, giving the
same lower bound for the distinct inert nodal domains meeting $\beta$.
\end{proof}

\begin{remark}We may also use nonnegative linear combinations of correlations. For finitely many $c_j\ge0$ and $\delta_j>0$, summing
\eqref{eq:negative-pair-bound} gives
\begin{equation}\label{eq:positive-combination-count}
 \sum_jc_j\int_{\RR}v(s)v(s+\delta_j)\,ds\le-D<0
 \quad\Longrightarrow\quad
 K\ge\frac{D}{S^2\sum_jc_j\delta_j}.
\end{equation}
For fixed $a_j>0$ and $c_j\ge0$ satisfying $\sum_jc_j\cos a_j=0$,
Theorem~\ref{prop:correlation} gives
\[
 \sum_jc_jI_\phi(a_j)
 =2H_\psi\sum_jc_jJ_0(a_j)+o(1).
\]
A negative limiting constant gives \eqref{eq:positive-combination-count}
with $\delta_j=\log(1+a_j/T)$.
\end{remark}

\section{Coefficient bounds and mean values}
\label{sec:analytic-inputs}

We collect the upper bounds used in the Mellin reduction and in the
ranges near $t=\pm T$. For each fixed $C>0$, the Fourier coefficient estimate
of Ghosh, Reznikov and Sarnak \cite[Prop.~3.7]{GRS13} gives
\begin{equation}
 \frac1{L(1,\sym^2\phi)}\sum_{n\le CT}|\lambda_\phi(n)|^2\ll_C T.
 \label{eq:normalised-reference-bound}
\end{equation}
Indeed, their coefficients are
$\widetilde\rho_\phi(n)=\rho_\phi(n)e^{-\pi T/2}$, and
\eqref{eq:fourier-normalization} gives
\[
 |\widetilde\rho_\phi(n)|^2
 =\frac{1+e^{-2\pi T}}{L(1,\sym^2\phi)}|\lambda_\phi(|n|)|^2.
\]
Their upper bound at the fixed length $CT$ therefore implies
\eqref{eq:normalised-reference-bound}. In particular,
\begin{equation}
 \frac1{\sqrt{L(1,\sym^2\phi)}}
 \sum_{n\le CT}\frac{|\lambda_\phi(n)|}{\sqrt n}
 \ll_C\sqrt{T\log T}
 \label{eq:mellin-coefficient-sum}
\end{equation}
by Cauchy--Schwarz. These bounds hold independently of the bulk limit
that we prove below.

The Montgomery--Vaughan mean-value theorem \cite{MV74}, in the form
\cite[Lem.~2.2]{YoungGeo}, implies that for any finite sequence
supported on $cN\le n\le CN$ and any interval $J$ of length at most $H$,
\begin{equation}
 \int_J\left|\sum_n a_n n^{-it}\right|^2dt
 \ll_C(H+N)\sum_n|a_n|^2.
 \label{eq:fixed-mean-value}
\end{equation}
The estimate is uniform in the location of $J$, since translating $t$
multiplies the coefficients by numbers of absolute value one. The
cutoffs arising here have the form $\psi(A(t)/n)$, with $\psi$ fixed.
We use Mellin inversion to separate the $t$- and $n$-dependence of this cutoff
and obtain the following adaptation.

\begin{lemma}\label{lem:fixed-profile-mean-value}
Fix $0<c<C$ and $\psi\in C_c^\infty((0,\infty))$. Let $N,H>0$,
let $J$ be an interval of length at most $H$, and let $(a_n)$ be
supported on $cN\le n\le CN$. For any measurable
$A:J\to(0,\infty)$,
\begin{equation}
 \int_J\left|\sum_n a_n n^{-it}\psi(A(t)/n)\right|^2dt
 \ll_{C,\psi}(H+N)\sum_n|a_n|^2.
 \label{eq:fixed-profile-mean-value}
\end{equation}
\end{lemma}

\begin{proof}
We put $\Psi(iv)=\int_0^\infty\psi(y)y^{iv}\,dy/y$.
By Mellin inversion, we obtain the exact separation
\[
 \sum_n a_n n^{-it}\psi(A(t)/n)
 =\frac1{2\pi}\int_\RR\Psi(iv)A(t)^{-iv}
       \sum_n a_n n^{-i(t-v)}\,dv.
\]
We apply Minkowski's inequality and \eqref{eq:fixed-mean-value} to the
inner polynomial on $J-v$. Since $|A(t)^{-iv}|=1$, the resulting
$L^2(J)$ bound is at most
\[
 \frac1{2\pi}\int_\RR|\Psi(iv)|\,dv\,
 \left(C_1(H+N)\sum_n|a_n|^2\right)^{1/2}.
\]
The integral is finite because $\psi$ is fixed, smooth and compactly
supported in $(0,\infty)$.
\end{proof}

We also use the compact Iwaniec--Sarnak bound \cite{IS95}: for every
fixed compact $\mathcal K\subset X$ and every $\eps>0$,
\begin{equation}\label{eq:supnorm}
 \sup_{\mathcal K}|\phi|\ll_{\mathcal K,\eps}T^{5/12+\eps}.
\end{equation}
In our norm-one convention this is
\cite[App.~C, Lem.~C.1, $k=0$]{GRS13}. We will use it to control the narrow
neighborhoods $\bigl||t|-T\bigr|\le T^\delta$ and convert a negative correlation into a count
of sign changes.

\section{The Mellin polynomial}
\label{sec:mellin}

We first express the correlation in \eqref{eq:arithmetic-correlation}
by the Mellin--Parseval formula. We then use the Fourier expansion and Mellin inversion to approximate
the Mellin transform of the weighted restriction by an explicit
gamma factor times a Dirichlet polynomial. Put $f(y)=\psi(y)\phi(iy)$ and use the convention
\[
 \mathcal M(s)=\int_0^\infty f(y)y^s\,\frac{dy}{y},\qquad
 \Psi(u)=\int_0^\infty\psi(y)y^u\,\frac{dy}{y}.
\]
For $g(s)=f(e^s)$, we have $\mathcal M(it)=\int_\RR g(s)e^{its}\,ds$.
Thus the shift $q=1+\alpha/T>0$ corresponds to translation by
$\log q$, and the Mellin--Parseval formula gives
\begin{equation}
 I_\phi(\alpha)=\frac1{2\pi}\int_\RR
 \cos(t\log q)|\mathcal M(it)|^2\,dt.
 \label{eq:mellin-correlation}
\end{equation}
Here $g$ is real, so $|\mathcal M(it)|^2$ is even and the sine term
integrates to zero. This identity includes both cutoffs in
\eqref{eq:arithmetic-correlation}.

To express $\mathcal M$ in terms of Dirichlet polynomials, we set
\[
 G_T(s)=\frac{\rho_\phi(1)}{2\pi^{s+1/2}}
 \Gamma\left(\frac{1/2+s+iT}{2}\right)
 \Gamma\left(\frac{1/2+s-iT}{2}\right).
\]
The Mellin integral \cite[10.43.19]{DLMF}
\[
 \int_0^\infty K_{iT}(v)v^{z-1}\,dv
 =2^{z-2}\Gamma((z+iT)/2)\Gamma((z-iT)/2),\qquad \Re z>0,
\]
combined with the Fourier expansion in \eqref{eq:fourier-normalization}, gives
\[
 \int_0^\infty\phi(iy)y^s\,\frac{dy}{y}=G_T(s)L(1/2+s,\phi)
\]
in an absolutely convergent half-plane. Applying Mellin inversion
to $\psi$, for fixed $c>2$, we obtain
\begin{equation}
 \mathcal M(it)=\frac1{2\pi i}\int_{\Re u=-c}
 \Psi(u)G_T(it-u)L(1/2+it-u,\phi)\,du.
 \label{eq:mellin-exact}
\end{equation}
On this line the Dirichlet series converges absolutely, the gamma
factors decay exponentially, and $\Psi$ decays faster than any power.
These facts justify the exchanges of summation and integration.

For $|t|<T$, the leading term of the gamma quotient localizes the
coefficients to $n\asymp\sqrt{T^2-t^2}$. We define
\begin{equation}
 P_T(t)=\sum_{n\ge1}\frac{\lambda_\phi(n)}{n^{1/2+it}}
       \psi\left(\frac{\sqrt{T^2-t^2}}{2\pi n}\right).
 \label{eq:mellin-polynomial}
\end{equation}
We need an approximation whose remainder can be controlled
in mean square: without an additional weight in the bulk, and with
the factor $O(H/T)$ on each boundary block near $t=\pm T$.
The same type of approximation of the squared Mellin transform under
the integral is used for Eisenstein series in
\cite[(2.11)--(2.12), Lem.~2.5]{YoungGeo}.
The estimates below use our cusp-form coefficients.

\begin{lemma}[Mellin approximation]\label{lem:mellin-reduction}
Fix $\delta>0$. Uniformly for $|t|\le T-T^\delta$,
\begin{equation}
 \mathcal M(it)=G_T(it)P_T(t)
 +O_{\psi,\delta}\left(
       T^{1/4}(T-|t|)^{-5/4}\sqrt{\log T}\right).
 \label{eq:mellin-finite}
\end{equation}
Moreover,
\begin{equation}
 |G_T(it)|^2=\frac{2\pi}{L(1,\sym^2\phi)\sqrt{T^2-t^2}}
       \{1+O((T-|t|)^{-1})\}.
 \label{eq:mellin-amplitude}
\end{equation}
\end{lemma}

\begin{proof}
We choose $a>0$ with $\supp\psi\subset[a,\infty)$. The integral
representation \cite[10.32.9]{DLMF} gives
$|K_{iT}(x)|\le K_0(x)\ll x^{-1/2}e^{-x}$ for $x\ge1$.
The $n=1$ term of \eqref{eq:normalised-reference-bound} gives
$|\rho_\phi(1)|e^{-\pi T/2}\ll\sqrt T$. Together with the crude
bound $|\lambda_\phi(n)|\ll n$, 
this permits us to truncate the Fourier sum defining $f$ at $|n|\le CT$,
with an exponentially small error in $\mathcal M(it)$, uniformly in
real $t$. Here we choose $C=C(\psi)$ so that $2\pi aC>\pi/2$,
ensuring exponential decay of the omitted tail. We also require
$C>1/(2\pi a)$ so that the truncation includes every nonzero term of
$P_T(t)$: by \eqref{eq:mellin-polynomial}, such a term has $n\le T/(2\pi a)$.

We write $h=T-|t|\ge T^\delta$ and $z_\pm=i(t\pm T)/2$.
Stirling's formula gives
\[
 \left|\Gamma\left(\frac{1}{4}+\frac{i(t+T)}{2}\right)
 \Gamma\left(\frac{1}{4}+\frac{i(t-T)}{2}\right)\right|^2
 =\frac{8\pi^2e^{-\pi T}}{\sqrt{T^2-t^2}}\{1+O(h^{-1})\}.
\]
The first-coefficient identity in \eqref{eq:fourier-normalization}
now proves \eqref{eq:mellin-amplitude}.
We let $\mathcal M_{\le CT}$ denote the Mellin transform of
\[
 2\rho_\phi(1)\psi(y)\sqrt y
 \sum_{1\le n\le CT}\lambda_\phi(n)K_{iT}(2\pi ny),
\]
where the factor $2$ combines the equal positive and negative Fourier
terms at $x=0$. By termwise Mellin inversion, we obtain the finite sum
\begin{equation}\label{eq:mellin-truncated-convolution}
 \mathcal M_{\le CT}(it)
 =\sum_{n\le CT}\frac{\lambda_\phi(n)}{n^{1/2+it}}
 \frac1{2\pi i}\int_{\Re u=-c}
       \Psi(u)G_T(it-u)n^u\,du.
\end{equation}
We shift each of these finitely many integrals to $\Re u=0$.
The decay of $\Psi$
controls the horizontal integrals. With $u=iv$ and $|v|\le h^{1/8}$,
logarithmic Stirling in a fixed sector
\cite[5.11.1, \S5.11(ii)]{DLMF} gives
\begin{align}
 &\frac{\Gamma(z_++1/4-iv/2)\Gamma(z_-+1/4-iv/2)}
 {\Gamma(z_++1/4)\Gamma(z_-+1/4)}=\left(\frac{\sqrt{T^2-t^2}}2\right)^{-iv}
       \left\{1+O\left(\frac{(1+|v|)^2}{h}\right)\right\}.
 \label{eq:gamma-uniform-remainder}
\end{align}
Its error integrates against $\Psi(iv)$ to $O_\psi(h^{-1})$.

For the complementary range of $v$, we use Stirling's bound
\cite[5.11.9]{DLMF}, extended over bounded $u$ by continuity:
\[
 |\Gamma(1/4+iu)|\asymp(1+|u|)^{-1/4}e^{-\pi|u|/2}
 \qquad(u\in\RR).
\]
It bounds the quotient in \eqref{eq:gamma-uniform-remainder}
by $O(T^{1/2})$.
Indeed its exponential factor is
\[
 \exp\left\{-\frac\pi4
 (|t+T-v|+|t-T-v|-2T)\right\}\le1,
\]
and its polynomial factor is at most $O((T^2-t^2)^{1/4})$.
The rapid decay of $\Psi(iv)$ therefore bounds this tail and the tail
of the leading term by $O_{\psi,K}(T^{1/2}h^{-K})$ for any fixed $K$.

Mellin inversion of the leading power gives the cutoff in $P_T(t)$.
By \eqref{eq:mellin-coefficient-sum} and \eqref{eq:mellin-amplitude},
the summed local error is at most
\[
 (Th)^{-1/4}h^{-1}\sqrt{T\log T}
       =T^{1/4}h^{-5/4}\sqrt{\log T}.
\]
The complementary tail contributes
$O_{\psi,K}(T^{3/4}h^{-K-1/4}\sqrt{\log T})$.
We choose $K>1+1/(2\delta)$ to absorb it into the same bound and prove
\eqref{eq:mellin-finite}.
\end{proof}

In the bulk calculation, we replace $|\mathcal M(it)|^2$ by its Dirichlet polynomial approximation. We therefore need
an integrated comparison of these squares on a fixed bulk interval. For
fixed $0<\eta<1$, the squared remainder in \eqref{eq:mellin-finite}
has integral $O_{\eta,\psi}(\log T/T)$ over $|t|\le(1-\eta)T$.
On this interval the union of the supports in $P_T$ lies in
$[c_{\eta,\psi}T,C_\psi T]$. We apply
Lemma~\ref{lem:fixed-profile-mean-value} with $N=T$ and
$A(t)=\sqrt{T^2-t^2}/(2\pi)$, truncating the coefficients to this
fixed interval. Equations \eqref{eq:mellin-amplitude} and
\eqref{eq:normalised-reference-bound} give
\[
 \int_{|t|\le(1-\eta)T}|G_T(it)P_T(t)|^2\,dt
 \ll_{\eta,\psi}\frac1{L(1,\sym^2\phi)}
       \sum_{cT\le n\le CT}\frac{|\lambda_\phi(n)|^2}{n}
 \ll_{\eta,\psi}1.
\]
We bound the cross term by Cauchy--Schwarz and use
\eqref{eq:mellin-amplitude} to replace $|G_T(it)|^2$ by its
leading term. The latter replacement contributes
$O_{\eta,\psi}(T^{-1})$ to the integral, and we obtain
\begin{align}
 &\int_{|t|\le(1-\eta)T}
 \left||\mathcal M(it)|^2-
       \frac{2\pi|P_T(t)|^2}
       {L(1,\sym^2\phi)\sqrt{T^2-t^2}}\right|\,dt\notag\\
 &\hspace{12mm}\ll_{\eta,\psi}\sqrt{\frac{\log T}{T}}.
 \label{eq:mellin-quadratic-error}
\end{align}
Near the endpoints we instead integrate the squared remainder with
the small multiplicative factor from the correlation difference, as in
\eqref{eq:boundary-mellin-error} below.

\section{The bulk contribution}
\label{sec:bulk}

On a fixed range $|t|\le(1-\eta)T$, with $\eta>0$, the polynomial
$P_T(t)$ has length comparable to $T$. We open its square to obtain shifted
convolution sums. For each fixed shift, the weight converges to a fixed
smooth function, and we will be able to apply arithmetic QUE.
The oscillation in the original kernel gives a summable bound for the
remaining shifts.

Lindenstrauss and Soundararajan proved full phase-space arithmetic QUE
\cite[Thm.~1.4]{Lin06}, \cite[Thm.~1]{Sound10}.
We use its qualitative shifted-convolution formulation recorded by Jung
\cite[(1.5) and Rem.~1.4]{JungSC}, which we state here in our normalization.

\begin{lemma}[Fixed-shift coefficient QUE]\label{lem:fixed-shift-que}
Let $h\in\ZZ$ be fixed, and let $W$ be a smooth function with compact
support in $(0,\infty)$, independent of $T$.
As $T\to\infty$ through the full family of forms in
Theorem~\ref{thm:main},
\begin{equation}\label{eq:fixed-shift-que}
 \frac1{T L(1,\sym^2\phi)}
 \sum_{\substack{n\ge1\\n+h\ge1}}
 \lambda_\phi(n)\lambda_\phi(n+h)W(n/T)
 \longrightarrow \frac6{\pi^2}\mathbf1_{h=0}
                         \int_0^\infty W(x)\,dx.
\end{equation}
For complex-valued $W$, we apply the real-valued statement to its real
and imaginary parts.
\end{lemma}
We prove Lemma~\ref{lem:fixed-shift-que} directly from full
phase-space QUE in Appendix~\ref{app:coefficient-que}.

To match the constant, Jung's first Fourier coefficient has square
$2/L(1,\sym^2\phi)$ by \eqref{eq:fourier-normalization}. We apply his
formula with weight $W(\,\cdot\,/\pi)$ and combine the positive and
negative indices using even parity. Its support stays away from zero,
so $n$ and $n+h$ have the same sign there for large $T$.
After reindexing, the two weights
are $W(n/T)$ and $W((n+h)/T)$. We replace $W((n+h)/T)$ by $W(n/T)$,
introducing an error of $O_{W,h}(T^{-1})$ in the normalized sum, by
\eqref{eq:normalised-reference-bound} and Cauchy--Schwarz.

The next proposition evaluates the correlation integral away from
$t=\pm T$. In the final application, we take $F(\tau)$ to be
$\cos(\alpha\tau)-\cos\alpha$ times a fixed even smooth cutoff
supported in $(-1,1)$. The boundary estimates then control
the omitted frequencies.

\begin{proposition}\label{prop:arithmetic-bulk}
For every fixed real even $F\in C_c^\infty((-1,1))$,
\begin{equation}\label{eq:bulk-target}
 \frac1{2\pi}\int_{\RR}F(t/T)|\mathcal M(it)|^2\,dt
 \longrightarrow
 \frac{2H_\psi}{\pi}
 \int_{-1}^1\frac{F(\tau)}{\sqrt{1-\tau^2}}\,d\tau.
\end{equation}
\end{proposition}

\begin{proof}
We choose fixed $\eta>0$ and $1<a<b$ so that
$\supp F\subset[-1+\eta,1-\eta]$ and $\supp\psi\subset[a,b]$.
Using the integrated comparison \eqref{eq:mellin-quadratic-error} and
changing variables by $t=T\tau$, we obtain
\begin{equation}\label{eq:bulk-total-reduction}
 \frac1{2\pi}\int_{\RR}F(t/T)|\mathcal M(it)|^2\,dt
 =\frac1{L(1,\sym^2\phi)}\int_{-1}^1
   \frac{F(\tau)|P_T(T\tau)|^2}{\sqrt{1-\tau^2}}\,d\tau+o_{F,\psi}(1).
\end{equation}
Opening the square we obtain $\sum_{h\in\ZZ}S_{T,h}$, where
\begin{align}
 S_{T,h}
 &=\frac1{L(1,\sym^2\phi)}
   \sum_{\substack{n\ge1\\n+h\ge1}}
   \frac{\lambda_\phi(n)\lambda_\phi(n+h)}{\sqrt{n(n+h)}}
                      \mathcal K_{T,F}(n,n+h),
 \label{eq:bulk-polynomial-shift}\\
 \mathcal K_{T,F}(n,m)
 &=\int_{-1}^1\frac{F(\tau)}{\sqrt{1-\tau^2}}
   \psi\!\left(\frac{T\sqrt{1-\tau^2}}{2\pi n}\right)
   \psi\!\left(\frac{T\sqrt{1-\tau^2}}{2\pi m}\right)
   e^{iT\tau\log(m/n)}\,d\tau.
 \label{eq:bulk-polynomial-kernel}
\end{align}
The cutoff factors restrict both $n$ and $m$ to $[cT,CT]$,
where $0<c<C$ depend only on $F$ and $\psi$.
We fix $h$ and put $x=n/T$. Uniformly on this compact range,
\[
 \frac1{\sqrt{x(x+h/T)}}=\frac1x+O_h(T^{-1}),\qquad
 T\log\left(1+\frac h{Tx}\right)=\frac hx+O_h(T^{-1}).
\]
We apply Taylor's formula to the second cutoff in
\eqref{eq:bulk-polynomial-kernel}, with $m=n+h$:
\[
 \psi\!\left(\frac{\sqrt{1-\tau^2}}{2\pi(x+h/T)}\right)
 =\psi\!\left(\frac{\sqrt{1-\tau^2}}{2\pi x}\right)
   +O_{F,\psi,h}(T^{-1}).
\]
The error is uniform for $x\in[c,C]$ and $\tau\in\supp F$.
Thus we replace the product of the two cutoffs by the square of the first.
Combining this with the preceding two approximations, we obtain the weight
\begin{equation}\label{eq:bulk-fixed-weight}
 W_{F,h}(x)=\frac1x\int_{-1}^1
  \frac{F(\tau)}{\sqrt{1-\tau^2}}
  \psi\!\left(\frac{\sqrt{1-\tau^2}}{2\pi x}\right)^2
  e^{ih\tau/x}\,d\tau.
\end{equation}
This weight is independent of $T$ and belongs to $C_c^\infty((0,\infty))$,
with support in $[c,C]$.
The error in replacing
$\mathcal K_{T,F}(n,n+h)/\sqrt{n(n+h)}$ by
$T^{-1}W_{F,h}(n/T)$ is $O_{F,\psi,h}(T^{-2})$.
This error remains small after summation: for a fixed enlarged $C'$,
\begin{align*}
 &\frac1{T L(1,\sym^2\phi)}
 \sum_{\substack{cT\le n\le CT\\n+h\ge1}}
       |\lambda_\phi(n)\lambda_\phi(n+h)|\\
 &\qquad\le \frac1{T L(1,\sym^2\phi)}
       \sum_{n\le C'T}|\lambda_\phi(n)|^2\ll_{F,\psi,h}1
\end{align*}
by \eqref{eq:normalised-reference-bound}. Consequently,
\begin{equation}\label{eq:bulk-frozen-shift}
 S_{T,h}=\frac1{T L(1,\sym^2\phi)}
 \sum_{\substack{n\ge1\\n+h\ge1}}
 \lambda_\phi(n)\lambda_\phi(n+h)W_{F,h}(n/T)
       +O_{F,\psi,h}(T^{-1}).
\end{equation}
We apply Lemma~\ref{lem:fixed-shift-que} with the fixed weight function
$W_{F,h}$. For every fixed $h\ne0$, the contribution $S_{T,h}$
tends to zero. For $h=0$, we substitute
$y=\sqrt{1-\tau^2}/(2\pi x)$ to obtain
\[
 \int_0^\infty W_{F,0}(x)\,dx
 =\frac{\pi H_\psi}{3}\int_{-1}^1\frac{F(\tau)}{\sqrt{1-\tau^2}}\,d\tau.
\]
The diagonal therefore has the limit in \eqref{eq:bulk-target}.

We still need to bound all large shifts together. In the original kernel
\eqref{eq:bulk-polynomial-kernel}, every fixed $\tau$-derivative of
the nonoscillatory factor is bounded uniformly on its support.
Integration by parts has no boundary terms, since $F$ is compactly
supported in $(-1,1)$. Moreover, $n,m\le CT$ gives
$T|\log(m/n)|\ge |m-n|/C$. Hence for every $A>0$,
\begin{equation}\label{eq:bulk-kernel-decay}
 |\mathcal K_{T,F}(n,m)|
 \ll_{A,F,\psi}(1+T|\log(m/n)|)^{-A}
 \ll_{A,F,\psi}(1+|m-n|)^{-A}.
\end{equation}
Using $\sqrt{nm}\ge cT$, Cauchy--Schwarz and
\eqref{eq:normalised-reference-bound}, we obtain
\[
 |S_{T,h}|\ll_{A,F,\psi}
 \frac{(1+|h|)^{-A}}{T L(1,\sym^2\phi)}
       \sum_{n\le CT}|\lambda_\phi(n)|^2
 \ll_{A,F,\psi}(1+|h|)^{-A}.
\]
The bound is uniform in $T$ and $h$.
In particular, for $A>1$,
\begin{equation}\label{eq:bulk-shift-tail}
 \sum_{|h|>K}|S_{T,h}|\ll_{A,F,\psi}K^{1-A},\qquad K\ge1.
\end{equation}
We may now keep $K$ fixed and let $T\to\infty$ in the finitely many
terms with $|h|\le K$. Only the diagonal survives. We then let
$K\to\infty$ in \eqref{eq:bulk-shift-tail} to prove the proposition.
\end{proof}

\section{The boundary range}
\label{sec:arithmetic-boundary}
On a boundary block $H\le T-|t|\le2H$, with $H\ge1$, the correlation difference
in \eqref{eq:arithmetic-defect} has an additional factor of size $O(H/T)$.
Indeed, writing
\begin{equation}\label{eq:defect-multiplier}
 b_{\alpha,T}(t)=\cos\bigl(t\log(1+\alpha/T)\bigr)-\cos\alpha
\end{equation}
we obtain for fixed $\alpha$ using Taylor's theorem
\begin{equation}\label{eq:boundary-taylor}
 |b_{\alpha,T}(t)|\ll_\alpha H/T+T^{-1}\ll_\alpha H/T.
\end{equation}

On a block of length $H$, the polynomial $P_T$ has length
$N=\sqrt{TH}$. We therefore need a bound for the normalized
coefficient sum at this shorter length. The following lemma
uses Soundararajan's comparison theorem \cite[Thm.~3]{Sound10}.

\begin{lemma}\label{lem:coefficient-bound}
Let $0<c<C$ be fixed and put $N=\sqrt{TH}$, where $1\le H\le T$.
For all sufficiently large $T$, uniformly in $H$,
\begin{equation}\label{eq:boundary-Q}
 Q_{c,C}(N):=\frac1{L(1,\sym^2\phi)}
       \sum_{cN\le n\le CN}\frac{|\lambda_\phi(n)|^2}{n}
 \ll_{c,C}(T/H)^{1/4}\log(eT/H).
\end{equation}
\end{lemma}

\begin{proof}
We put $S(X)=\sum_{n\le X}|\lambda_\phi(n)|^2$.
The Hecke relation and $\lambda_\phi(1)=1$ allow us to apply
Soundararajan's comparison theorem \cite[Thm.~3]{Sound10}:
\[
 S(x/y)\ll (1+\log y)y^{-1/2}S(x),\qquad 1\le y\le x.
\]
We take $x=CT$ and $y=T/N$, so that $x/y=CN$.
Since $\sqrt T\le N\le T$, we have $1\le y\le x$ for all
sufficiently large $T$, uniformly in $H$. Equation
\eqref{eq:normalised-reference-bound} at the fixed length $CT$ gives
\begin{equation}\label{eq:boundary-coefficients}
 \frac{S(CN)}{L(1,\sym^2\phi)}
 \ll_C\sqrt{TN}\log(eT/N)
 \ll_C T^{3/4}H^{1/4}\log(eT/H).
\end{equation}
Dividing by $cN$, we obtain \eqref{eq:boundary-Q}. The reference multiple
$C$ is independent of $H$, so the threshold is uniform over the blocks.
\end{proof}
We now bound the boundary contribution uniformly in $\eta$.
\begin{proposition}\label{prop:arithmetic-boundary}
Fix $\alpha\in\RR$ and $0<\delta<1/6$. Uniformly for
$0<\eta\le1/4$ and all sufficiently large $T$,
\begin{equation}\label{eq:arithmetic-boundary}
 \frac1{2\pi}
 \int_{T^\delta\le T-|t|\le\eta T}
 |b_{\alpha,T}(t)|\,|\mathcal M(it)|^2\,dt
 \ll
 \eta^{3/4}\log(e/\eta)+T^{-(1+\delta)/2}\log T.
\end{equation}
The implied constant and the threshold for $T$ depend only on
$\alpha$, $\psi$ and $\delta$.
\end{proposition}

\begin{proof}
We consider $J_H=[T-2H,T-H]$, or a subinterval, where
$T^\delta\le H\le\eta T$, and put $N=\sqrt{TH}$.
Since $H\le\eta T\le T/4$, on $J_H$ we have
\[
 T^2-t^2=(T-t)(T+t)\asymp HT,
 \qquad \frac{\sqrt{T^2-t^2}}{2\pi}\asymp N,
\]
with absolute implied constants. On $J_H$, the support of $\psi$
restricts the sum to $c_\psi N\le n\le C_\psi N$,
where $0<c_\psi<C_\psi$ are fixed. We apply
Lemma~\ref{lem:fixed-profile-mean-value} to
$a_n=\lambda_\phi(n)n^{-1/2}$ restricted to this interval,
with $A(t)=\sqrt{T^2-t^2}/(2\pi)$, to obtain
\[
 \int_{J_H}|P_T(t)|^2\,dt
 \ll_\psi (H+N)L(1,\sym^2\phi)Q_{c_\psi,C_\psi}(N).
\]
Combining this estimate with \eqref{eq:mellin-amplitude} and
\eqref{eq:boundary-taylor}, we obtain
\begin{align}
 \int_{J_H}|b_{\alpha,T}(t)|\,|G_T(it)P_T(t)|^2\,dt
 &\ll_{\alpha,\psi}
 \frac HT\,\frac{H+N}{N}\,Q_{c_\psi,C_\psi}(N)
 \notag\\
 &\ll_{\alpha,\psi}(H/T)^{3/4}\log(eT/H).
 \label{eq:arithmetic-boundary-block}
\end{align}
Here $H/N=\sqrt{H/T}\le1$. The factor $H/T$ in \eqref{eq:boundary-taylor}
thus absorbs the coefficient loss $(T/H)^{1/4}$ in
\eqref{eq:boundary-Q}.

Next we estimate the error $\mathcal M(it)-G_T(it)P_T(t)$
with the weight $|b_{\alpha,T}(t)|$.
By Lemma~\ref{lem:mellin-reduction}, its squared absolute value is
$O_{\psi,\delta}(T^{1/2}H^{-5/2}\log T)$ on $J_H$.
Multiplying by $H/T$ and integrating over an interval of length
at most $H$, we obtain
\begin{equation}\label{eq:boundary-mellin-error}
 \int_{J_H}|b_{\alpha,T}(t)|
       |\mathcal M(it)-G_T(it)P_T(t)|^2\,dt
 \ll_{\alpha,\psi,\delta}\frac{\log T}{\sqrt{TH}}.
\end{equation}
The inequality $|u+v|^2\le2|u|^2+2|v|^2$ bounds the weighted
integral of $|\mathcal M(it)|^2$ on $J_H$ by twice the sum of
the left-hand sides of \eqref{eq:arithmetic-boundary-block}
and \eqref{eq:boundary-mellin-error}.

We cover $[T^\delta,\eta T]$ by the dyadic intervals $[H,2H]$,
$H=2^jT^\delta$, truncating the last one. Using the bound \eqref{eq:arithmetic-boundary-block} we get
\begin{align}
 \sum_{\substack{H=2^jT^\delta\le\eta T\\j\ge0}}
       (H/T)^{3/4}\log(eT/H)
 &\ll\eta^{3/4}\sum_{j\ge0}2^{-3j/4}
                  \bigl(\log(e/\eta)+j\bigr)
 \notag\\
 &\ll\eta^{3/4}\log(e/\eta).
 \label{eq:boundary-dyadic-sum}
\end{align}
Summing \eqref{eq:boundary-mellin-error} over the same blocks gives
$O_{\alpha,\psi,\delta}(T^{-(1+\delta)/2}\log T)$,
uniformly in $\eta$. Since $b_{\alpha,T}$ and $|\mathcal M(it)|^2$ are even,
the negative-frequency blocks have the same bound.
\end{proof}

\section{Completion of the correlation formula}
\label{sec:arithmetic-completion}

We fix $0<\delta<1/6$. The bulk and boundary estimates leave the
two-sided neighborhoods \mbox{$\bigl||t|-T\bigr|\le T^\delta$} and the
exterior $|t|>T+T^\delta$. In the neighborhoods, the multiplicative factor
$b_\alpha(t/T)$ defined below is small enough to use the compact supremum bound.
In the exterior, the gamma factors and the cutoff transform give
rapid decay.

\begin{lemma}\label{lem:exterior-tail}
For every fixed $A>0$ and $r\ge0$,
\[
 \int_{|t|>T+T^\delta}(1+|t|)^r|\mathcal M(it)|^2\,dt
 \ll T^{-A}.
\]
The implied constant depends only on $A$, $r$, $\delta$ and $\psi$.
\end{lemma}
\begin{proof}
We use \eqref{eq:mellin-exact} with $c=3$.
The elementary Hecke bound gives $|L(7/2+iu,\phi)|\ll1$ uniformly in
$\phi$ and $u$, while \eqref{eq:normalised-reference-bound} and
\eqref{eq:fourier-normalization} give
$|\rho_\phi(1)|e^{-\pi T/2}\ll T^{1/2}$.
By Stirling's bound \cite[5.11.9]{DLMF}, we obtain
\begin{equation}\label{eq:global-mellin-envelope}
 |G_T(3+iu)L(7/2+iu,\phi)|
 \ll(1+T+|u|)^3
       \exp\left(-\frac\pi2(|u|-T)_+\right),
\end{equation}
where $(x)_+=\max(x,0)$.

For $d=|t|-T>0$, we write $u=-3+iv$ in \eqref{eq:mellin-exact}
and split the integral at $|v|=d/2$.
The first part has the exponential factor $e^{-\pi d/4}$ from
\eqref{eq:global-mellin-envelope}. In the second, we use the rapid
decay of $\Psi(-3+iv)$ and
$1+T+|t-v|\le(1+T+|t|)(1+|v|)$. For every fixed $B>0$ this gives
\[
 |\mathcal M(it)|\ll_{B,\psi}
 (1+T+|t|)^3(1+|t|-T)^{-B},\qquad |t|>T.
\]
We square this bound and integrate separately over
$T^\delta<d\le T$ and $d>T$. We conclude by taking $B$ sufficiently
large in terms of $A,r,\delta$.
\end{proof}
\begin{lemma}\label{lem:core-and-phase}
For fixed $\alpha\in\RR$, put
$b_\alpha(\tau)=\cos(\alpha\tau)-\cos\alpha$.
Then
\begin{align*}
 \int_{\bigl||t|-T\bigr|\le T^\delta}
 |b_\alpha(t/T)|\,|\mathcal M(it)|^2\,dt&=o(1),\\
 \int_\RR |b_{\alpha,T}(t)-b_\alpha(t/T)|
                    |\mathcal M(it)|^2\,dt&=o(1).
\end{align*}
\end{lemma}

\begin{proof}
The compact supremum estimate \eqref{eq:supnorm} and Parseval give,
for every $\varepsilon>0$,
\begin{equation}\label{eq:restriction-mass-upper}
 I_\phi(0)=\frac1{2\pi}\int_\RR|\mathcal M(it)|^2\,dt
 \ll_{\psi,\varepsilon}T^{5/6+\varepsilon}.
\end{equation}
The stronger unconditional bound
$I_\phi(0)\ll_{\psi,\varepsilon}T^\varepsilon$ was proven in \cite{GRS13}.
The elementary estimate \eqref{eq:restriction-mass-upper} suffices here.
On $\bigl||t|-T\bigr|\le T^\delta$, we apply Taylor's theorem to obtain
$|b_\alpha(t/T)|\ll_\alpha T^{\delta-1}$.
The first integral is therefore
$O_{\alpha,\psi,\varepsilon}(T^{\delta-1/6+\varepsilon})=o(1)$
provided we choose $0<\varepsilon<1/6-\delta$.

We use $\log(1+\alpha/T)=\alpha/T+O_\alpha(T^{-2})$ to obtain
$|b_{\alpha,T}(t)-b_\alpha(t/T)|\ll_\alpha T^{-1}$ for $|t|\le2T$.
By \eqref{eq:restriction-mass-upper}, this part of the second integral is
$O_{\alpha,\psi,\varepsilon}(T^{-1/6+\varepsilon})=o(1)$.
On $|t|>2T$, we use
Lemma~\ref{lem:exterior-tail} instead.
\end{proof}

\begin{proof}[Proof of Theorem~\ref{prop:correlation}]
By the Mellin identity \eqref{eq:mellin-correlation} and
Lemma~\ref{lem:core-and-phase}, we have
\[
 I_\phi(\alpha)-\cos\alpha\,I_\phi(0)
 =\frac1{2\pi}\int_\RR b_\alpha(t/T)|\mathcal M(it)|^2\,dt+o(1).
\]
We fix $0<\eta<1/8$ and choose an even
$w_\eta\in C_c^\infty((-1,1))$, with $0\le w_\eta\le1$, equal
to one for $|\tau|\le1-2\eta$ and zero for $|\tau|\ge1-\eta$.
We evaluate the part with $w_\eta(t/T)$ by
Proposition~\ref{prop:arithmetic-bulk}, with the fixed test function
$F(\tau)=w_\eta(\tau)b_\alpha(\tau)$.

For the complementary part, we divide the frequencies into
\[
 T^\delta<T-|t|<2\eta T,\qquad
 \bigl||t|-T\bigr|\le T^\delta,\qquad
 |t|>T+T^\delta.
\]
On the first range, Proposition~\ref{prop:arithmetic-boundary}
applies with $2\eta\le1/4$. Since $0\le1-w_\eta\le1$,
the same absolute bound holds after inserting $1-w_\eta(t/T)$.
Lemma~\ref{lem:core-and-phase} allows us to replace
$b_{\alpha,T}(t)$ by $b_\alpha(t/T)$ with an $o(1)$ error
independent of $\eta$.
The same lemma controls the second range, including the exterior
halves of the neighborhoods. Lemma~\ref{lem:exterior-tail}
controls the third exterior range. Hence
\begin{align}\label{eq:arithmetic-ordered-limit}
 I_\phi(\alpha)-\cos\alpha\,I_\phi(0)
 ={}&\frac{2H_\psi}{\pi}\int_{-1}^1
       \frac{w_\eta(\tau)b_\alpha(\tau)}{\sqrt{1-\tau^2}}\,d\tau
       +o_{\alpha,\psi,\eta}(1)\notag\\
 &+O_{\alpha,\psi}\bigl(\eta^{3/4}\log(e/\eta)\bigr).
\end{align}
The last implied constant is independent of $T$ and $\eta$.
In the bulk calculation, for each fixed $\eta$ we first fix finitely
many coefficient shifts, take $T\to\infty$, and then remove
the shift cutoff using its uniform summable tail. We may now let
$\eta\downarrow0$.

The boundary error vanishes in this limit. Since
$b_\alpha(\tau)=O_\alpha(1-|\tau|)$ at both endpoints, we remove
$w_\eta$ from the limiting integral with an error of $O_\alpha(\eta^{3/2})$.
Finally,
\[
 \frac1\pi\int_{-1}^1
       \frac{b_\alpha(\tau)}{\sqrt{1-\tau^2}}\,d\tau
       =J_0(\alpha)-\cos\alpha,
\]
which proves \eqref{eq:arithmetic-defect}. At $\alpha=0$ this
identity is tautological, so it gives no asymptotic for $I_\phi(0)$.
\end{proof}

\appendix
\section{The fixed-shift formulation of QUE}
\label{app:coefficient-que}

We derive Lemma~\ref{lem:fixed-shift-que} from full phase-space
arithmetic QUE. Throughout this appendix we fix the shift $h$ and every
weight before $T\to\infty$.

We write the Fourier expansion as
\[
 \phi(x+iy)=\sum_{n\in\ZZ}c_{T,n}(y)e(nx),\qquad
 c_{T,n}(y)=\rho_\phi(1)\lambda_\phi(|n|)\sqrt y
                  K_{iT}(2\pi|n|y),\quad c_{T,0}=0.
\]
We first choose a fixed bounded open interval $J$ whose closure lies in
$(0,\infty)$, such
that the cylinder $(\RR/\ZZ)\times\overline J$ contains no lift
of an elliptic point. Such intervals exist arbitrarily close to
zero: the heights of those lifts form a locally finite subset of
$(0,\infty)$. The image of the closed cylinder in $X$ is compact. For
$b\in C_c^\infty(J)$ and every fixed integer $j\ge0$, full QUE gives
\begin{align}
 &\int_J b(y)\sum_n(n/T)^j c_{T,n}(y)
                   \overline{c_{T,n+h}(y)}\,\frac{dy}{y^2}
 \notag\\
 &\qquad\longrightarrow
 \frac3\pi\mathbf1_{h=0}\int_J b(y)\frac{dy}{y^2}
       \frac1{2\pi}\int_0^{2\pi}
       \left(\frac{\cos\theta}{2\pi y}\right)^j d\theta.
 \label{eq:que-cylinder-moments}
\end{align}
We justify this use of QUE by localization. The left side is
\[
 \int_J\int_0^1 b(y)e(hx)
 \left(\frac{\partial_x}{2\pi iT}\right)^j
 \phi(x+iy)\overline{\phi(x+iy)}\,dx\,\frac{dy}{y^2}.
\]
We use a smooth partition of unity to localize the integral to
finitely many patches on which the projection to $X$ is injective.
On each patch, the displayed differential operator, with its
coefficient extended by zero, is compactly supported in a smooth
chart of $X$. An input cutoff equal to one near the support of that coefficient
leaves the expression unchanged. The leading symbol is the
corresponding partition weight times
$b(y)e(hx)(\xi_x/(2\pi))^j$.
Fixed local elliptic estimates bound every normalized derivative
of $\phi$ in $L^2$ on these charts. We choose the input cutoff with
support in a larger chart neighborhood. Its commutator with
$T^{-2}\Delta_+-(1+1/(4T^2))$ is supported away from the output
coefficient. A fixed frequency cutoff equal to one near
$y^2(\xi_x^2+\xi_y^2)=1$ reduces the operator to order zero.
On the complementary frequencies a local elliptic parametrix
applies, and the separated commutator support makes its contribution
$O_N(T^{-N})$ by the same local Sobolev bounds. Lower-order terms
from composition are $O(T^{-1})$.
This is the local finite-order extension of QUE described in
\cite{JJ18}. It applies on these compact charts with
$T^{-2}(T^2+1/4)=1+O(T^{-2})$.
We apply microlocal QUE to the resulting order-zero operators
and sum over the partition.
The Liouville integral unfolds back to the cylinder, counting each
piece once even if their images overlap. On the unit cotangent
circle $\xi_x=\cos\theta/y$, and
$\int_0^1e(hx)\,dx=\mathbf1_{h=0}$, proving
\eqref{eq:que-cylinder-moments}.

The differential moments above determine limits for polynomial
weights in $n/T$. To pass to arbitrary fixed smooth weights,
we fix $0\le\chi_0\in C_c^\infty(J)$ and define measures on
$J\times\RR$ by
\[
 d\nu_{T,h}(y,s)=\chi_0(y)\sum_n
 c_{T,n}(y)\overline{c_{T,n+h}(y)}
 \delta_{n/T}(ds)\,\frac{dy}{y^2}.
\]
The measure $\nu_{T,0}$ is positive. If $\tau_a$ translates the
second coordinate by $a$, then
\[
 |\nu_{T,h}|\le\tfrac12\bigl(\nu_{T,0}
                             +\tau_{-h/T}\nu_{T,0}\bigr).
\]
Equation~\eqref{eq:que-cylinder-moments} bounds all even moments
of $\nu_{T,0}$ and gives their limits. Thus these measures and
$\nu_{T,h}$ are tight and have bounded total variation.
If $a=\inf\supp\chi_0>0$, their subsequential limits are supported
in $|s|\le(2\pi a)^{-1}$. Indeed, for any $R>(2\pi a)^{-1}$,
the $2m$-th moment bounds the limiting mass outside $[-R,R]$
by $C(2\pi aR)^{-2m}$, which tends to zero as $m\to\infty$.
The domination above gives the same conclusion for $\nu_{T,h}$.
Higher even moments also make each fixed polynomial moment
uniformly integrable, so its limit passes to every subsequential
limit. Polynomial approximation on the resulting compact set,
together with smooth weights in $y$, determines that limit.
Consequently, for every $A\in C_c^\infty(J\times\RR)$, we choose
$\chi_0=1$ on its $y$-support and obtain
\begin{align}
 &\int_J\sum_n A(y,n/T)c_{T,n}(y)
           \overline{c_{T,n+h}(y)}\,\frac{dy}{y^2}
 \notag\\
 &\qquad\longrightarrow
 \frac3\pi\mathbf1_{h=0}\int_J\frac{dy}{y^2}
   \frac1{2\pi}\int_0^{2\pi}
   A\left(y,\frac{\cos\theta}{2\pi y}\right)d\theta.
 \label{eq:que-cylinder-weights}
\end{align}
We fix all moment orders and approximation weights before taking the
spectral limit. The cutoff $\chi_0$ equals one on the required
support, so its insertion changes neither side.

Next we choose a weight that extracts $W(n/T)$ from the paired
Bessel integrals. Their leading product contains a cosine factor;
we arrange for its weighted $y$-integral to recover $W$.
We now fix $W\in C_c^\infty((0,\infty))$, with
$\supp W\subset[c,C]$ for $0<c<C$. We choose the interval $J$ above
small enough that $2\pi C\sup J<1$, and take
$0\le\chi\in C_c^\infty(J)$, $\chi\not\equiv0$. On a fixed
neighborhood of $\supp W$ we put
\[
 r(y,s)=\sqrt{1-(2\pi sy)^2},\qquad
 B_h(s)=\int_J\chi(y)\cos^2\!\left(\frac{h r(y,s)}s\right)
                                  \frac{dy}{y}.
\]
The function $B_h$ is smooth and strictly positive. For $h=0$ this is
immediate. For $h\ne0$, the argument $h r(y,s)/s$ has nonzero
$y$-derivative, so its cosine cannot vanish on an interval where
$\chi>0$. In particular $B_h$ has a positive minimum on
$\supp W$. The fixed smooth weight
\[
 A(y,s)=\frac{\chi(y)W(s)r(y,s)}{B_h(s)}
                    \cos\!\left(\frac{h r(y,s)}s\right),
\]
extended by zero, belongs to $C_c^\infty(J\times(0,\infty))$.

For completeness, we record the oscillatory Bessel asymptotic used here:
\[
 e^{\pi T/2}K_{iT}(Tz)
 =\sqrt{\frac{2\pi}{T}}(1-z^2)^{-1/4}
   \sin\!\left(T\left(\operatorname{arcosh}(1/z)
                    -\sqrt{1-z^2}\right)+\frac\pi4\right)
   +O(T^{-3/2}),
\]
uniformly for $z$ in a fixed compact subinterval of $(0,1)$
\cite[Lem.~3.1]{GRS13}.
We set $s=n/T$ in the support of $W$. The two Bessel phases at
$2\pi sy$ and $2\pi(s+h/T)y$ differ by
$-h r(y,s)/s+O(T^{-1})$. Their sum has $y$-derivative
$-T(r(y,s)+r(y,s+h/T))/y$, whose absolute value is bounded
below by a fixed positive multiple of $T$. Multiplying the
asymptotics and integrating once by parts in the sum phase,
we therefore obtain, uniformly in these $n$,
\begin{align*}
 &\int_J A(y,s)K_{iT}(2\pi ny)K_{iT}(2\pi(n+h)y)
                                  \frac{dy}{y}\\
 &\quad=\frac{\pi e^{-\pi T}}T
   \int_J\frac{A(y,s)}{r(y,s)}
                   \cos\!\left(\frac{h r(y,s)}s\right)\frac{dy}{y}
          +O(e^{-\pi T}T^{-2})\\
 &\quad=\frac{\pi e^{-\pi T}}T W(s)
          +O(e^{-\pi T}T^{-2}).
\end{align*}
There are no boundary terms because $\chi$ is compactly supported.
The error combines the uniform Bessel remainder, the fixed-shift
Taylor error and the integration-by-parts error.
We use \eqref{eq:fourier-normalization} to evaluate the leading
coefficient, and \eqref{eq:normalised-reference-bound} together
with Cauchy--Schwarz to sum the error terms. This gives
\begin{align}
 &\int_J\sum_n A(y,n/T)c_{T,n}(y)
          \overline{c_{T,n+h}(y)}\,\frac{dy}{y^2}
 \notag\\
 &\qquad=
 \frac\pi{T L(1,\sym^2\phi)}
   \sum_{\substack{n\ge1\\n+h\ge1}}
   \lambda_\phi(n)\lambda_\phi(n+h)W(n/T)+O(T^{-1}).
 \label{eq:que-extract-coefficients}
\end{align}
Indeed, $\rho_\phi(1)^2e^{-\pi T}
=(1+e^{-2\pi T})/L(1,\sym^2\phi)$, and the normalized absolute
coefficient-pair sum over $n\asymp T$ is bounded.

For $h\ne0$, \eqref{eq:que-cylinder-weights} makes the left side
of \eqref{eq:que-extract-coefficients} tend to zero. For $h=0$,
we substitute $s=\cos\theta/(2\pi y)$ on the two intervals
where $\cos\theta>0$ to obtain its limit as
\[
 \frac6\pi\int_J\int_0^\infty
        \frac{A(y,s)}{r(y,s)}\,ds\,\frac{dy}{y}
 =\frac6\pi\int_0^\infty W(s)\,ds,
\]
since $B_0=\int_J\chi(y)\,dy/y$. Dividing by $\pi$ in
\eqref{eq:que-extract-coefficients}, we obtain
\eqref{eq:fixed-shift-que}, including its constant.

\providecommand{\bysame}{\leavevmode\hbox to3em{\hrulefill}\thinspace}
\providecommand{\MR}{\relax\ifhmode\unskip\space\fi MR }
\providecommand{\MRhref}[2]{\href{http://www.ams.org/mathscinet-getitem?mr=#1}{#2}
}
\providecommand{\href}[2]{#2}

\end{document}